\documentclass[12pt,a4paper,twoside]{article}

\newcommand{\pageformat}[6]{\setlength{\hoffset}{-1in}
                  \setlength{\voffset}{-1in}
                  \addtolength{\hoffset}{#5}
                            \addtolength{\voffset}{#6}
                            \setlength{\oddsidemargin}{#1}
                            \setlength{\evensidemargin}{#2}
                            \setlength{\textwidth}{\paperwidth}
                  \addtolength{\textwidth}{-\oddsidemargin}
                  \addtolength{\textwidth}{-\evensidemargin}
                  \addtolength{\textwidth}{-\marginparsep}
                  \addtolength{\textwidth}{-\marginparwidth}
                            \setlength{\topmargin}{#3}
                            \setlength{\textheight}{\paperheight}
                  \addtolength{\textheight}{-\topmargin}
                  \addtolength{\textheight}{-\headheight}
                  \addtolength{\textheight}{-\headsep}
                  \addtolength{\textheight}{-\footskip}
                  \addtolength{\textheight}{-#4}}
\pageformat{2cm}{3cm}{25mm}{25mm}{1pt}{0pt}

\usepackage{ifthen}
\newboolean{article}
    \setboolean{article}{true}
\newboolean{report}
\newboolean{book}
\newboolean{letter}
\newboolean{german}
\newboolean{italian}
\newboolean{nobaselinestretch}
\newboolean{nosectionappendix}
\newboolean{oldtoc}
\newboolean{nosectionequn}
\newboolean{notheorem}

\ifthenelse{\boolean{german}}{
    \usepackage{german}}{}

\usepackage[latin1]{inputenc}

\usepackage{amsmath}
\usepackage{amssymb}
\usepackage[mathscr]{eucal}

\ifthenelse{\boolean{notheorem}}{}{
    \usepackage{theorem}}

\ifthenelse{\boolean{nobaselinestretch}}{}{
    \renewcommand{\baselinestretch}{1.25}}

\newenvironment{env}[2]{\begin{#1}#2\end{#1}}{}
    \newcommand{\beq}[1]{\begin{env}{equation}{#1}}
    \newcommand{\beqn}[1]{\begin{env}{equation*}{#1}}
    \newcommand{\bal}[1]{\begin{env}{align}{#1}}
    \newcommand{\baln}[1]{\begin{env}{align*}{#1}}
    \newcommand{\bga}[1]{\begin{env}{gather}{#1}}
    \newcommand{\bgan}[1]{\begin{env}{gather*}{#1}}
    \newcommand{\bflal}[1]{\begin{env}{flalign}{#1}}
    \newcommand{\bflaln}[1]{\begin{env}{flalign*}{#1}}
    \newcommand{\bmu}[1]{\begin{env}{multline}{#1}}
    \newcommand{\bmun}[1]{\begin{env}{multline*}{#1}}
    \newcommand{\bsp}[1]{\begin{env}{split}{#1}}

    \newcommand{\eeq}{\end{env}}
    \newcommand{\eeqn}{\end{env}}
    \newcommand{\eal}{\end{env}}
    \newcommand{\ealn}{\end{env}}
    \newcommand{\ega}{\end{env}}
    \newcommand{\egan}{\end{env}}
    \newcommand{\eflal}{\end{env}}
    \newcommand{\eflaln}{\end{env}}
    \newcommand{\emu}{\end{env}}
    \newcommand{\emun}{\end{env}}
    \newcommand{\esp}{\end{env}}

\newcommand{\lf}{\vspace{2ex}}

\renewcommand{\bf}[1]{\textbf{#1}}
\renewcommand{\it}[1]{\textit{#1}}

\renewcommand{\sf}[1]{\textsf{#1}}

\renewcommand{\tt}[1]{\texttt{#1}}
\newcommand{\hl}[1]{\bf{\it{#1}}}

\newcommand{\mbf}[1]{\mathbf{#1}}
\newcommand{\msf}[1]{\text{\small$\sf{#1}$}}

\newcommand{\cmc}[1]{\mathcal{#1}}
\newcommand{\eus}[1]{\mathscr{#1}}
\newcommand{\euf}[1]{\mathfrak{#1}}
\newcommand{\bb}[1]{\mathbb{#1}}

\newcommand{\mscriptsize}[1]{{\setlength{\arraycolsep}{.3ex}\text{\scriptsize$#1$}}}

\newcommand{\nbd}[1]{$#1$\nobreakdash--}
\newcommand{\ol}[1]{\overline{#1}}

\newcommand{\wh}[1]{\widehat{#1}}
\newcommand{\ve}{\varepsilon}
\newcommand{\vt}{\vartheta}

\newcommand{\vp}{\varphi}

\newcommand{\norm}[1]{\left\lVert#1\right\rVert}

\newcommand{\snorm}[1]{\norm{\smash{#1}}}

\newcommand{\family}[1]{\left(#1\right)}

\newcommand{\bfam}[1]{\bigl(#1\bigr)}

\newcommand{\AB}[1]{\langle#1\rangle}

\newcommand{\CB}[1]{\{#1\}}
\newcommand{\bCB}[1]{\bigl\{#1\bigr\}}
\newcommand{\BCB}[1]{\Bigl\{#1\Bigr\}}
\newcommand{\SB}[1]{[#1]}

\newcommand{\RO}[1]{[#1)}

\newcommand{\Matrix}[1]{\begin{pmatrix}#1\end{pmatrix}}

\newcommand{\sMatrix}[1]{\mscriptsize{\Matrix{#1}}}

\newcommand{\set}[2][]{
    \ifthenelse{\equal{#1}{}}{
        \CB{#2}}{
        \CB{#1~|~#2}}}
\newcommand{\bset}[2][]{
    \ifthenelse{\equal{#1}{}}{
        \bCB{#2}}{
        \bCB{#1~|~#2}}}
\newcommand{\Bset}[2][]{
    \ifthenelse{\equal{#1}{}}{
        \BCB{#2}}{
        \BCB{#1~\big|~#2}}}

\DeclareMathOperator{\ls}{\normalfont\msf{span}}
\DeclareMathOperator{\cls}{\ol{\ls}}

\DeclareMathOperator{\id}{\normalfont\msf{id}}

\newcommand{\N}{\bb{N}}

\newcommand{\R}{\bb{R}}
\newcommand{\bS}{\bb{S}}

\newcommand{\cB}{\cmc{B}}

\newcommand{\sB}{\eus{B}}

\newcommand{\sK}{\eus{K}}

\newcommand{\sN}{\eus{N}}

\newcommand{\sS}{\eus{S}}

\newcommand{\en}{\euf{n}}

\newcommand{\eR}{\euf{R}}

\newcommand{\U}{\mbf{1}}

\ifthenelse{\boolean{nosectionequn}}{}{
    \numberwithin{equation}{section}
    }

\ifthenelse{\boolean{article}\or\boolean{letter}\or\boolean{nosectionequn}}{
    \setboolean{nosectionappendix}{true}}{}
\ifthenelse{\boolean{nosectionappendix}}{}{
    \renewcommand{\appendix}{
        \chapter*{\appendixname}
        \addcontentsline{toc}{chapter}{\appendixname}
        \renewcommand{\thesection}{\Alph{section}}
        \setcounter{section}{0}}}
   
\ifthenelse{\boolean{report}\or\boolean{book}}{
    }{}

\ifthenelse{\boolean{notheorem}}{}{
        \newcommand{\mnname}{Mathematical note.}
        \newcommand{\enname}{End of the note.}
        \newcommand{\definame}{Definition.}
        \newcommand{\propname}{Proposition.}
        \newcommand{\lemname}{Lemma.}
        \newcommand{\exname}{Example.}
        \newcommand{\exername}{Exercise.}
        \newcommand{\remname}{Remark.}
        \newcommand{\obname}{Observation.}
        \newcommand{\thmname}{Theorem.}
        \newcommand{\corname}{Corollary.}
        \newcommand{\proofname}{Proof.}
    \ifthenelse{\boolean{german}}{
        \renewcommand{\mnname}{Mathematische Notiz.}
        \renewcommand{\enname}{Ende der Notiz.}
        \renewcommand{\exname}{Beispiel.}
        \renewcommand{\exername}{Übung.}
        \renewcommand{\remname}{Bemerkung.}
        \renewcommand{\obname}{Beobachtung.}
        \renewcommand{\thmname}{Satz.}
        \renewcommand{\corname}{Korollar.}
        \renewcommand{\proofname}{Beweis.}}{}
    \ifthenelse{\boolean{italian}}{
        \renewcommand{\mnname}{Nota matematica.}
        \renewcommand{\enname}{Fina della nota.}
        \renewcommand{\definame}{Definizione.}
        \renewcommand{\propname}{Proposizione.}
        \renewcommand{\exname}{Esempio.}
        \renewcommand{\exername}{Esercizio.}
        \renewcommand{\remname}{Nota.}
        \renewcommand{\obname}{Osservazione.}
        \renewcommand{\thmname}{Teorema.}
        \renewcommand{\corname}{Corollario.}
        \renewcommand{\proofname}{Dimostrazione.}

       \renewcommand{\appendixname}{Appendice}

       }{}
    \theoremheaderfont{\normalfont\bfseries}
    \theoremstyle{change}
        \theorembodyfont{\rmfamily}
            \newtheorem{emp}{}[section]
                \newcommand{\bemp}[1][]{
                    \begin{emp}\hskip-\labelsep\bf{#1}\hskip\labelsep}
                \newcommand{\eemp}{\end{emp}}
\newtheorem{itemp}[emp]{}
                \newcommand{\bitemp}[1][]{
                    \begin{itemp}\hskip-\labelsep\bf{#1}\hskip\labelsep\normalfont\itshape}
                \newcommand{\eitemp}{\end{itemp}}
            \newtheorem{mn}[emp]{\mnname}
                \newcommand{\bnm}{\begin{mn}~\begin{quotation}\renewcommand{\baselinestretch}{1}\small\noindent\ignorespaces}
                \newcommand{\enm}{\end{quotation}\hfill\bf{\enname}\end{mn}}
            \newtheorem{ex}[emp]{\exname}
                \newcommand{\bex}{\begin{ex}}
                \newcommand{\eex}{\end{ex}}
            \newtheorem{exer}[emp]{\exername}
                \newcommand{\bexer}{\begin{exer}}
                \newcommand{\eexer}{\end{exer}}
            \newtheorem{defi}[emp]{\definame}
                \newcommand{\bdefi}{\begin{defi}}
                \newcommand{\edefi}{\end{defi}}
            \newtheorem{rem}[emp]{\remname}
                \newcommand{\brem}{\begin{rem}}
                \newcommand{\erem}{\end{rem}}
            \newtheorem{ob}[emp]{\obname}
                \newcommand{\bob}{\begin{ob}}
                \newcommand{\eob}{\end{ob}}
        \theorembodyfont{\normalfont\itshape}
            \newtheorem{thm}[emp]{\thmname}
                \newcommand{\bthm}{\begin{thm}}
                \newcommand{\ethm}{\end{thm}}
            \newtheorem{prop}[emp]{\propname}
                \newcommand{\bprop}{\begin{prop}}
                \newcommand{\eprop}{\end{prop}}
            \newtheorem{cor}[emp]{\corname}
                \newcommand{\bcor}{\begin{cor}}
                \newcommand{\ecor}{\end{cor}}
            \newtheorem{lem}[emp]{\lemname}
                \newcommand{\blem}{\begin{lem}}
                \newcommand{\elem}{\end{lem}}
\newenvironment{empn}[1]{\lf\noindent\bf{#1}\ignorespaces\hskip\labelsep}{\lf}
		\newcommand{\bempn}[1]{\begin{empn}{#1}}
		\newcommand{\eempn}{\end{empn}}
		\newcommand{\bitempn}[1]{\begin{empn}{#1}\normalfont\itshape}
		\newcommand{\eitempn}{\end{empn}}
                \newcommand{\bnmn}{\begin{empn}{\mnname}~\begin{quotation}\renewcommand{\baselinestretch}{1}\small\noindent\ignorespaces}
                \newcommand{\enmn}{\end{quotation}\hfill\bf{\enname}\end{empn}}
		\newcommand{\bexn}{\begin{empn}{\exname}}
		\newcommand{\eexn}{\end{empn}}
		\newcommand{\bexern}{\begin{empn}{\exername}}
		\newcommand{\eexern}{\end{empn}}
		\newcommand{\bdefin}{\begin{empn}{\definame}}
		\newcommand{\edefin}{\end{empn}}
		\newcommand{\bremn}{\begin{empn}{\remname}}
		\newcommand{\eremn}{\end{empn}}
		\newcommand{\bobn}{\begin{empn}{\obname}}
		\newcommand{\eobn}{\end{empn}}

\newcommand{\qedsymbol}{~\rule[-0.35mm]{2mm}{2mm}}
    \newcounter{proof}[emp]
    \newenvironment{Proof}[1]{
        \vspace{1ex}
        \renewcommand{\item}[1][\stepcounter{proof}(\roman{proof})]%
            {##1\hskip\labelsep}
        \noindent\textsc{#1\hskip\labelsep}}{
        \nolinebreak\qedsymbol}
    \newcommand{\proof}[1][\proofname]{
        \begin{Proof}{#1}\ignorespaces}
    \newcommand{\qed}{\end{Proof}}
    \newcommand{\noqed}{
        \renewcommand{\qedsymbol}{}
        \end{Proof}}}
    \ifthenelse{\boolean{italian}}{
        \renewcommand{\proofname}{Dimostrazione.}}{}

\usepackage[varg]{txfonts}

\usepackage[hypertex]{hyperref}

\begin{document}

\bibliographystyle{amsalpha}

\title{$E_0$--Semigroups for Continuous Poduct Systems:\\The Nonunital Case}
\author{}
\author{
~\\
Michael Skeide\thanks{This work is supported by research funds of University of Molise and Italian MIUR.}\\\\
{\small\itshape Universit\`a\ degli Studi del Molise}\\
{\small\itshape Dipartimento S.E.G.e S.}\\
{\small\itshape Via de Sanctis}\\
{\small\itshape 86100 Campobasso, Italy}\\
{\small{\itshape E-mail: \tt{skeide@unimol.it}}}\\
{\small{\itshape Homepage: \tt{http://www.math.tu-cottbus.de/INSTITUT/lswas/\_skeide.html}}}\\
}
\date{January 2009}

{
\renewcommand{\baselinestretch}{1}
\maketitle


\begin{abstract}
\noindent
Let $\cB$ be a \nbd{\sigma}unital \nbd{C^*}algebra. We show that every strongly continuous \nbd{E_0}semigroup on the algebra of adjointable operators on a full Hilbert \nbd{\cB}module $E$ gives rise to a full continuous product system of correspondences over $\cB$. We show that every full continuous product system of correspondences over $\cB$ arises in that way. If the product system is countably generated, then $E$ can be chosen countable generated, and if $E$ is countably generated, then so is the product system. We show that under these countability hypotheses there is a one-to-one correspondence between \nbd{E_0}semigroup up to stable cocycle conjugacy and continuous product systems up isomorphism. This generalizes the results for unital $\cB$ to the \nbd{\sigma}unital case.
\end{abstract}

}



\section{Introduction}\label{intro}

Factorizable families of Hilbert spaces are known since quite a while; see Araki \cite{Ara70}, Streater \cite{Str69}, and Parthasarathy and Schmidt \cite{PaSchm72}. Arveson \cite{Arv89,Arv90a,Arv89a,Arv90} developed this idea into a concise theory of \it{tensor product systems} of Hilbert spaces (\hl{Arveson systems}, for short). Roughly speaking, Arveson's theory provides a classification of \nbd{E_0}semigroups (unital endomorphism semigroups) on $\sB(H)$ ($H$ a Hilbert space) by Arveson systems up to \it{cocycle conjugacy}. It is comparably plain to associate with every \nbd{E_0}semigroup an Arveson system, and to show that two \nbd{E_0}semigroups have isomorphic Arveson systems, if and only if they are cocycle conjugate. All this and an index theory for \nbd{E_0}semigroups is done in \cite{Arv89}. To show that every Arveson system comes from an \nbd{E_0}semigroup, was done the remaining articles \cite{Arv90a,Arv89a,Arv90}.

Liebscher \cite{Lie00p1} provided the second proof of this \it{fundamental theorem about Arveson systems}. This proof is still quite involved. But it adds the information that the \nbd{E_0}semigroup having the given Arveson system may be chosen \it{pure}. Only recently,  \cite{Ske06} we provided a simple and self-contained proof. Shortly later, Arveson \cite{Arv06} presented another simple proof. In \cite{Ske06a} we showed that the output of \cite{Arv06} and (a special case) of \cite{Ske06} are unitarily equivalent.

Meanwhile, several authors investigated tensor product systems of Hilbert bimodules or \it{correspondences}; see Bhat and Skeide \cite{BhSk00}, Muhly and Solel \cite{MuSo02}, and Hirschberg and Zacharias \cite{Hir04,HiZa03}. A connection between \nbd{E_0}semigroups on $\sB^a(E)$, the algebra of all adjointable operators on a full Hilbert \nbd{\cB}module $E$, and product systems of correspondences over $\cB$ (paralleling that of Arveson) has been established in Skeide \cite{Ske02} and, in its general version, in Skeide \cite{Ske04p'}.

Our scope that resulted in the simple proof of \cite{Ske06}, was to find a proof that works also for Hilbert modules. Funnily enough, in the two cases we could treat so far, namely when $\cB$ is a unital \nbd{C^*}algebra \cite{Ske07}, or when $\cB$ is a von Neumann algebra \cite{Ske07p} (in preparation), we proceeded utilizing Arveson's idea \cite{Arv06} in an essential way. In these notes, we now add the nonunital case under countability assumptions. ($\cB$ should be \nbd{\sigma}unital. And for the complete classification result, the occurring modules should be countably generated.) For the discrete case in Section \ref{discalgSEC} we need the original idea of \cite{Ske06}. The correct adaptation of Arveson's idea is a new crucial ingredient for the continuous time case. We also mention that the proof for unital $\cB$ in \cite{Ske07} that the \nbd{E_0}semigroup constructed there from a continuous product system induces the same continuous structure on that product system, contains a gap. The new Theorem \ref{cisothm} is far more general and fixes also the gap in \cite{Ske07}.

\section{The product system of an $E_0$--semigroup}\label{E0PSSEC}

Let $\bS$ denote one of the semigroups $\N_0=\CB{0,1,\ldots}$ and $\R_+=\RO{0,\infty}$. Fix a Hilbert \nbd{\cB}module $E$. In Skeide \cite{Ske02} we constructed the product system of a strict \nbd{E_0}semigroup $\vt$ on $\sB^a(E)$, following Bhat's construction in \cite{Bha96}, starting from a \hl{unit vector} $\xi\in E$, that is, from vector with ``length'' $\AB{\xi,\xi}=\U\in\cB$. This means, in particular, that $\cB$ is unital and that $E$ is full. The construction in \cite{Ske02} goes as follows.

Put $E_t:=\vt_t(\xi\xi^*)E$. Turn it into a correspondence over $\cB$ by defining the left action $bx_t:=\vt_t(\xi b\xi^*)x_t$. (Note that this left action is unital.) Define a map $v_t\colon E\odot E_t\rightarrow E$ by setting $v_t(x\odot y_t):=\vt_t(x\xi^*)y_t$. It is easy to check that this map is \hl{isometric} (that is, inner product preserving) and, therefore, well-defined. Surjectivity follows from strictness; see \cite{Ske02} for details. One easily verifies the following properties.
\begin{enumerate}
\item\label{Eps1}
$\vt$ can be recovered from the unitaries $v_t$ as $\vt_t(a)=v_t(a\odot\id_t)v_t^*$.

\item\label{Eps2}
The restriction $u_{s,t}$ to $E_s\odot E_t\subset E\odot E_t$ defines a bilinear unitary onto $E_{s+t}\subset E$.

\item\label{Eps3}
$E_0=\cB$ and $v_0$, $u_{s,0}$, and $u_{0,t}$ are the canonical identifications, that is, right multiplication (in the case of $v_0$ and $u_{s,0}$) and left multiplication (in the case of $u_{0,t}$), with the elements in $E_0=\cB$.

\item\label{Eps4}
Both ``multiplications'' $(x,y_t)\mapsto xy_t:=v_t(x\odot y_t)$ and $(x_s,y_t)\mapsto x_sy_t:=u_{s,t}(x_s\odot y_t)$ iterate associatively, that is, $(xy_s)z_t=x(y_sz_t)$ and $(x_ry_s)z_t=x_r(y_sz_t)$.
\end{enumerate}

\noindent
A family $E^\odot=\bfam{E_t}_{t\in\bS}$ of \nbd{\cB}correspondences with \hl{structure maps} $u_{s,t}$ fulfilling \ref{Eps2} and the relevant part of \ref{Eps3} and \ref{Eps4}, has been called \hl{product system} in Bhat and Skeide \cite{BhSk00}. Given a product system $E^\odot$ with structure maps $u_{s,t}$, a \hl{full} Hilbert \nbd{\cB}module $E$ (that is, the \hl{range ideal} $\cB_E:=\cls\AB{E,E}$ of $E$ coincides with $\cB$) and a family of unitaries $v_t$ fulfilling the relevant part of \ref{Eps3} and \ref{Eps4}, has been called a \hl{left dilation} of $E^\odot$ to $E$ in Skeide \cite{Ske07}. Note that if there exists a left dilation of $E^\odot$, then $E^\odot$ is necessarily \hl{full}, that is, $E_t$ is full for every $t$. If the $v_t$ form a left dilation of $E^\odot$ to $E$, then $\vt^v_t(a):=v_t(a\odot\id_t)v_t^*$ defines an \nbd{E_0}semigroup $\vt^v$. If $E$ has a unit vector, then the product system constructed from $\vt^v$ is (isomorphic to) $E^\odot$. Recall that a \hl{morphism} between product systems $E^\odot$ and $F^\odot$ is a family $w^\odot=\bfam{w_t}_{t\in\bS}$ of bilinear adjointable maps $w_t\colon E_t\rightarrow F_t$ such that $(w_sx_s)(w_ty_t)=w_{s+t}(x_sy_t)$ and $w_0=\id_\cB$. An \hl{isomorphism} is a morphism that consists of unitaries.

We say a strict \nbd{E_0}semigroup $\vt$ and a product system $E^\odot$ are \hl{associated}, if there exists a left dilation $v_t$ of $E^\odot$ such that $\vt=\vt^v$. It is known that for each strict \nbd{E_0}semigroup there is, up to isomorphism, only one product system that can be associated with that \nbd{E_0}semigroup; see Skeide \cite[Section 6]{Ske08p1}.

We have just seen that every strict \nbd{E_0}semigroup can be associated with a product system, provided that $E$ has a unit vector. There is a general construction in Skeide \cite{Ske04p'} for arbitrary (full) $E$ even if $\cB$ is nonunital, based on the representation theory of $\sB^a(E)$ from Muhly, Skeide, and Solel \cite{MSS06}. For the converse result, we have several stages:

\begin{enumerate}
\item
If $\cB$ is unital, we have the existence result \cite[Theorem 7.6]{Ske04p'} for the discrete case $\bS=\N_0$.

\item
Without continuity conditions, we can prove the continuous time case $\bS=\R_+$ by the method invented in Skeide \cite{Ske06} for the Hilbert space case, by reducing it to preceding result for the discrete case. Since, in the noncontinuous case, there are involved direct sums over the index set $\RO{0,1}$ and the shift on that set, the constructed \nbd{E_0}semigroup is definitely noncontinuous.

\item
In \cite{Ske07} we resolved, still for unital $\cB$, the continuous time case with continuity conditions both on the \nbd{E_0}semigroup and on the product system.

\end{enumerate}
(In Skeide \cite{Ske07p} we deal with the general von Neumann case. But this is out of the scope of the present notes, where we restrict to the \nbd{C^*}case.)

We see that in all three stages the case of nonunital $\cB$ is still missing. As for all three stages it is crucial to find a good adaptation of Bhat's method of the construction of the product system from an \nbd{E_0}semigroup (and not the abstract one based on \cite{MSS06}), we spend the present section to such find such a construction. In the following sections we apply the new insight to adapt also the proof for three stages, well, not to the general nonunital case, but to the \nbd{\sigma}unital case.

\lf
The crucial observation which gives the correct hint and resolves all problems that, so far, prevented us from dealing with the case of nonunital $\cB$, is quite simple. What does it mean if $E$ has a unit vector $\xi$? Well, it means that $E$ has a direct summand $\xi\cB\cong\cB$. The projection onto that summand is $\xi\xi^*$. If $E$ is full over a unital \nbd{C^*}algebra, then a finite multiple $E^n$ of $E$ will have a unit vector; see \cite[Lemma 3.2]{Ske04p'}. This was enough to treat the problems in the unital case. Now Lance \cite[Proposition 7.4]{Lan95} asserts the following: If $E$ is a full Hilbert module over a \nbd{\sigma}unital \nbd{C^*}algebra $\cB$, then $E^\infty$ has, well, not a unit vector, but a direct summand $\cB$. And this turns out to be enough for all our purposes.

To begin with let $\cB$ be an arbitrary \nbd{C^*}algebra. Suppose $E$ has a direct summand $\cB$, that is, suppose $E=\cB\oplus F$. This decomposes $\sB^a(E)$ into $\sMatrix{\sB^a(\cB)&\sB^a(F,\cB)\\\sB^a(\cB,F)&\sB^a(F)}$. Let $p\in\sB^a(E)$ denote the projection $((\beta,y)\mapsto(\beta,0)$ onto $\cB\subset E$. For $x\in E$ we define the element $xp\in\sB^a(E)$ by setting $xp(\beta,y):=x\beta$. The adjoint map is $px^*\colon x'\mapsto(\AB{x',x},0)$. Observe that $x'ppx*$ is just the usual rank-one operator $x'x^*$. Note, too, that $\pi\colon\AB{x,x'}\mapsto px^*x'p$ defines nothing but the canonical embedding of $\cB$ into the \nbd{\sB^a(\cB)}corner of $\sB^a(E)$.

Let $\vt$ be a strict \nbd{E_0}semigroup on $\sB^a(E)$. Following the procedure in presence of a unit vector, we put $E_t:=\vt_t(p)E$. It follows that $bx_t:=\vt_t(\pi(b))x_t$ defines a nondegenerate ($\vt_t$ is strict!) left action of $\cB$ on $E_t$ turning, thus, $E_t$ into a correspondence over $\cB$. By
\beqn{
v_t(x\odot y_t)
~:=~
\vt_t(xp)y_t
}\eeqn
we define a unitary $E\odot E_t\rightarrow E$. (By $\AB{\vt_t(xp)y_t,\vt_t(x'p)y'_t}=\AB{y_t,\AB{x,x'}y'_t}=\AB{x\odot y_t,x'\odot y'_t}$ we see that $v_t$ is isometric. Surjectivity follows from $\vt_t(xy^*)z=\vt_t(xp)\vt_t(py^*)z=v_t(x\odot(\vt_t(py^*)z)$, existence of a bounded approximate unit of finite-rank operators for $\sK(E)$ and strictness of $\vt_t$, in precisely the same way as in \cite{Ske02}.) Obviously, $\vt_t(a)=v_t(a\odot\id_t)v_t^*$. (Simply, apply both sides to $v_t(x\odot y_t)$.) The restriction $u_{s,t}$ of $v_t$ to $E_s\odot E_t$ is surjective onto $E_{s+t}$. (It is into, because $\vt_{s+t}(p)v_t(x\odot y_t)=v_t(\vt_s(p)x\odot y_t)$. It is onto, because $(\U-\vt_{s+t}(p))v_t(x\odot y_t)=v_t((\U-\vt_s(p))x\odot y_t)$.) Also the marginal conditions for $t=0$ or $s=0$ are satisfied. So, $E^\odot=\bfam{E_t}_{t\in\bS}$ is a product system and the $v_t$ form a left dilation giving back $\vt$ as $\vt^v$.

If $E$ has no direct summand but $\cB$ is \nbd{\sigma}unital, then we know that $E^\infty$ has a direct summand $\cB$. It is known that $\vt$ and its amplification $\vt^\infty$ to $\sB^a(E^\infty)$ have the same product system; see \cite[Section 9]{Ske08p1}. We, thus, proved the following.

\bthm
Let $\vt$ be a strict \nbd{E_0}semigroup on $\sB^a(E)$ where $E$ is a full Hilbert module over a \nbd{\sigma}unital \nbd{C^*}algebra. Then the product system of $\vt$ can be obtained by the prescribed construction applied to the amplification of $\vt$ to $\sB^a(E^\infty)$ based on any choice of a direct summand $\cB$ of $E^\infty$.
\ethm

If, in the continuous time case, $\vt$ is strongly continuous, then we would like that this property is reflected by a continuous structure of the product system. In Skeide \cite{Ske03b,Ske07} a \hl{continuous product system} is defined as a product system $E^\odot=\bfam{E_t}_{t\in\R_+}$ with a family of isometric embeddings $i_t\colon E_t\rightarrow\wh{E}$ into a right Hilbert \nbd{\cB}module $\wh{E}$ (there is no left action on $\wh{E}$) fulfilling the following conditions: Denote by
\beqn{
CS_i(E^\odot)
~=~
\BCB{\,\bfam{x_t}_{t\in\R_+}\colon x_t\in E_t,~t\mapsto i_tx_t\text{~is continuous~}}
}\eeqn
the set of \hl{continuous sections} of $E^\odot$ (with respect to the embeddings $i_t$). Then, firstly,
\beqn{
\BCB{\,x_s\colon\bfam{x_t}_{t\in\R_+}\in CS_i(E^\odot)\,}
~=~
E_s
}\eeqn
for all $s\in\R_+$ (that is, $E^\odot$ has sufficiently many continuous sections), and, secondly,
\beqn{
(s,t)
~\longmapsto~
i_{s+t}(x_sy_t)
}\eeqn
is continuous for all $\bfam{x_t}_{t\in\R_+},\bfam{y_t}_{t\in\R_+}\in CS_i(E^\odot)$ (that is, the `product' of continuous sections is continuous). A morphism between continuous product systems is \hl{continuous}, if it sends continuous sections to continuous sections. An \hl{isomorphism} of continuous product systems is a continuous isomorphism with continuous inverse. Clearly, an isomorphism provides a bijection between the sets of continuous sections.

The following theorem also settles a gap in the proof of \cite[Proposition 4.9]{Ske07} and generalizes it considerably. We illustrate its applications in the end of Section \ref{contSEC}.

\bthm\label{cisothm}
Let $i_t\colon E_t\rightarrow E^i$ and $k_t\colon E_t\rightarrow E^k$ be two continuous structures on the product system $E^\odot=\bfam{E_t}_{\in\R_+}$. If the identity morphism is a continuous morphism from $E^\odot$ with respect to the embeddings $i$ to $E^\odot$ with respect to the embeddings $k$, then the identity morphism is a continuous isomorphism.
\ethm

\proof
This statement means that if $x\in CS_i(E^\odot)$ $\Longrightarrow$ $x\in CS_k(E^\odot)$, then $x\in CS_i(E^\odot)$ $\Longleftrightarrow$ $x\in CS_k(E^\odot)$. Note that this is only a statement on the Banach bundle structure of $E^\odot$, while the product system structure does not play any role. Notice also the notion of uniform convergence of a sequence of sections of any subset $I$ of $\R_+$ depends only on the pointwise norms on $E_t$. It does not refer in any way to the embeddings $i_t$ or $k_t$. Nevertheless, a uniform limit on $I$ of sections that are continuous with respect to $i$ (to $k$) is continuous on $I$ with respect to $i$ (to $k$). Therefore, if we can approximate a section $x\in CS_k(E^\odot)$ on each compact interval $I=\SB{a,b}\subset\R_+$ uniformly by sections in $CS_i(E\odot)$, then $x\in CS_i(E^\odot)$.

Suppose $x$ and $I$ are as stated. For every $\beta\in I$ choose a section $y^\beta\in CS_i(E^\odot)\subset CS_k(E^\odot)$ such that $y^\beta_\beta=x_\beta$. Choose $\ve>0$. For every $\beta\in I$ choose an interval $I_\beta\subset I$ which is open in $I$ and which contains $\beta$ such that $\snorm{x_\alpha-y^\beta_\alpha}<\ve$ for all $\alpha\in I_\beta$. (Since $\snorm{x_\alpha-y^\beta_\alpha}=\snorm{k_\alpha x_\alpha-k_\alpha y^\beta_\alpha}$ and since $y^\beta,x\in CS_k(E^\odot)$, such $I_\beta$ exist.) So, we may choose $\beta_1,\ldots,\beta_m$ such that the union over $I_{\beta_i}$ is $\SB{a,b}$. By standard theorems about \it{partitions of unity} there exist continuous functions $\vp_i$ on $\SB{a,b}$ with the following properties:
\baln{
0
~\le~
&
\vp_i
~\le~
1,
&
\vp_i\upharpoonright I_{\beta_i}^\complement
&
~=~
0,
&
\sum_{i=1}^m\vp_i
&
~=~
1.
}\ealn
From these properties one verifies easily that $\norm{x_\alpha-\sum_{i=1}^m\vp_i(\alpha)y^{\beta_i}_\alpha}<\ve$ for all $\alpha\in I$. Since $\sum_{i=1}^m\vp_i(\alpha)y^{\beta_i}_\alpha\in CS_i(E^\odot)$, the section $x^\odot$ is the locally uniform limit of section in $CS_i(E^\odot)$ and, therefore, in $CS_i(E^\odot)$ itself.\qed

\lf
Already in \cite{Ske03b} we have shown that the product system $E^\odot$ of a strongly continuous strict \nbd{E_0}semigroup on $\sB^a(E)$ when $E$ has a unit vector, can be equipped with a continuous structure in the following way: Put $\wh{E}:=E$. It is, then, easy to see that $E^\odot$ is a continuous product system with respect to the canonical embeddings $i_t$ of the submodules $E_t\subset E$ into $E$. We now, simply, do the same for the product system constructed above in the case of general $\cB$ (and full $E$, of course).

Construct the amplification $\vt^\infty$ of $\vt$ on $\sB^a(\infty)$, so that $E^\infty$ has now a direct summand $\cB$ with projection $p\in\sB^a(E^\infty)$ onto that summand. Put $E_t:=\vt^\infty_t(p)E^\infty$ and choose for $i_t$ the canonical embeddings of $E_t$ into $E^\infty$. Precisely as in \cite{Ske08p1} (where the unital case has been treated, so that $p=\xi\xi^*$ for some unit vector $\xi\in E^\infty$) one shows that $E^\odot$ is continuous, that the continuous structure does not depend on the choice of the summand $\cB$ in $E^\infty$, and that, if $E$ has already a direct summand $\cB$, then the continuous structure is the same as if we had proceeded without amplifying $\vt$. We do not repeat the proof from \cite{Ske08p1} as it generalizes word by word.

This concludes the description of the construction of full product systems from strict $E_0$--semi\-groups $\vt=\bfam{\vt_t}_{t\in\bS}$ and of continuous full product systems from strongly continuous strict \nbd{E_0}semigroups $\vt=\bfam{\vt_t}_{t\in\R_+}$. The remainder of these notes is dedicated to the reverse direction.

\section{Discrete case and algebraic continuous time case}\label{discalgSEC}

Let $F$ be a full correspondence over $\cB$. We seek a full Hilbert \nbd{\cB}module $E$ such that $E\cong E\odot F$, for in that case this induces a unital strict endomorphism $\theta\colon a\mapsto a\odot\id_F\in\sB^a(E\odot F)\cong\sB^a(E)$ of $\sB^a(E)$ and the discrete product system $E^\odot$ associated with the discrete \nbd{E_0}semigroup $\vt=\bfam{\vt_n}_{n\in\N_0}$ with $\vt_n:=\theta^n$ is $E^\odot=\bfam{E_n}_{n\in\N_0}$ with $E_n:=F^{\odot n}$ and the canonical identifications $F^{\odot m}\odot F^{\odot n}=F^{\odot (m+n)}$.

If $E_1=F$ has a unit vector $\xi_1=\zeta$, then $\xi^\odot=\bfam{\xi_n}_{n\in\N_0}$ with $\xi_n:=\zeta^{\odot n}$ is a \it{unital unit}. In general, a \hl{unit} $\xi^\odot$ for a product system is a family $\xi^\odot$ of elements $\xi_t\in E_t$ such that $\xi_s\xi_t=\xi_{s+t}$ $(s,t\in\bS)$ and $\xi_0=\U\in\cB=E_0$. A unit is \hl{unital}, if it consists of unit vectors. Already Arveson \cite{Arv89} noted that in presence of a unital unit in a product system it is easy to construct an \nbd{E_0}semigroup associated with that product system. Simply, embed $E_t$ as $\xi_sE_t$ into $E_{s+t}$. These embeddings form an  inductive system and the factorization $u_{s,t}\colon E_s\odot E_t\rightarrow E_{s+t}$ ``survive'' the inductive limit as $v_t\colon E_\infty\odot E_t\rightarrow E_{\text{``$\infty+t$''}}=E_\infty$. Clearly, all associativity conditions are preserved so that the $v_t$ define a left dilation of $E^\odot$ to $E_\infty$ and the induced \nbd{E_0}semigroup $\vt^v$ has $E^\odot$ as product system.

The basic idea in Skeide \cite{Ske04p'} was: Even if $F$ has no unit vector, then $F^n$ has one for suitably big $n\in\N$. The same is true \it{cum grano salis} for the correspondence $M_\infty(F^n)\cong M_{n\cdot\infty,\infty}(F)\cong M_{\infty,\infty}(F)=M_\infty(F)$ over $M_\infty(\cB)$. The \it{cum grano salis} refers to that $M_\infty(\cB)=\sK(\cB^\infty)$ is always nonunital, and $M_{\infty}(F)=\sK(\cB^\infty,F^\infty)$ cannot contain a unit vector. What we need are both \it{strict completions} the multiplier algebra $\sB^a(\cB^\infty)$ of $\sK(\cB^\infty)$ and the correspondence $\sB^a(\cB^\infty,F^\infty)$ over $\sB^a(\cB^\infty)$. This correspondence \bf{does} have a unit vector. We find an inductive limit $E_\infty$ (a Hilbert \nbd{\sB^a(\cB^\infty)}module!) and a strict \nbd{E_0}semigroup $\vt^\infty$. If we define the (full) Hilbert \nbd{\cB}module $E:=E_\infty\odot\cB^\infty$, the \nbd{E_0}semigroup $\vt^\infty$ gives rise to a strict \nbd{E_0}semigroup $\vt$ on $\sB^a(E)\cong\sB^a(E_\infty)$, and $\vt$ has $E^\odot$ as product system.

All this has been described in \cite[Section 7]{Ske04p'} in a very detailed manner for unital $\cB$. The point is now that everything goes through precisely as in \cite{Ske04p'}, if $\cB$ is \nbd{\sigma}unital. Just that now instead of $F^n$ for $n\in\N$ we have to start with $F^\infty$. (\cite[Proposition 7.4]{Lan95} guarantees that $F^\infty$ has a direct summand $\cB$ and from this it follows that the correspondence $\sB^a(\cB^\infty,F^\infty)$ over $\sB^a(\cB^\infty)$ has a unit vector.) The rest goes precisely as in \cite{Ske04p'}. We, thus, proved:

\bthm\label{discthm}
Let $E^\odot=\bfam{E_n}_{n\in\N_0}$ be a full product system of correspondences over a \nbd{\sigma}unital \nbd{C^*}algebra. Then there exist a full Hilbert \nbd{\cB}module $E$ and a strict \nbd{E_0}semigroup $\vt$ on $\sB^a(E)$ such that the product system of $\vt$ is $E^\odot$.
\ethm

Once the discrete case (Theorem \ref{discthm}) is known, we may use it to construct a solution for the algebraic (without continuity conditions) continuous time case $\bS=\R_+$. This can be done by the procedure invented in Skeide \cite{Ske06} for the Hilbert space case. And we pointed it out in \cite{Ske07} for modules over unital $\cB$. The idea is the following: To find a left dilation of a full product system $E^\odot=\bfam{E_t}_{t\in\R_+}$ we start with a left dilation of the discrete subsystem $\bfam{E_t}_{t\in\N_0}$ to a Hilbert module $\breve{E}$, that is, with a family of unitaries $\breve{v}_n\colon\breve{E}\odot E_n\rightarrow\breve{E}$ that fulfill the necessary associativity conditions. We put $E:=\breve{E}\odot\int_0^1 E_\alpha\,d\alpha$. If $\cB$ is \nbd{\sigma}unital, then existence of such a dilation is granted by Theorem \ref{discthm}. The following identifications
\bal{
\notag
E\odot E_t
~=~
\breve{E}\odot\family{\int_0^1E_\alpha\,d\alpha}\odot E_t
&
~=~
\breve{E}\odot\int_t^{1+t}E_\alpha\,d\alpha
\\[2ex]
\notag
&
~\cong~
\family{\breve{E}\odot E_n\odot\int_{t-n}^1E_\alpha\,d\alpha}
\oplus
\family{\breve{E}\odot E_{n+1}\odot\int_0^{t-n}E_\alpha\,d\alpha}
\\[2ex]\label{idea}
&
~\cong~
\family{\breve{E}\odot\int_{t-n}^1E_\alpha\,d\alpha}
\oplus
\family{\breve{E}\odot\int_0^{t-n}E_\alpha\,d\alpha}
~=~
E
}\eal
suggest, then, a family of unitaries $v_t\colon E\odot E_t\rightarrow E$. The slightly tedious thing in \cite{Ske06} was to show associativity, that is, to show that the $v_t$ form a dilation to $E$.

By that method, whenever we are able to dilate the discrete subsystem $\bfam{E_t}_{t\in\N_0}$ of $E^\odot$ and to give a meaning to the direct integrals with respect to a translation (mod $1$) invariant measure, we are also able to dilate the whole product system $E^\odot$. In absence of continuity conditions, this translation invariant measure can only be the counting measure, so that the direct integrals are simply direct sums. We find:

\bthm\label{algcthm}
Let $E^\odot=\bfam{E_t}_{t\in\R_+}$ be a full product system of correspondences over a \nbd{\sigma}unital \nbd{C^*}algebra. Then there exist a full Hilbert \nbd{\cB}module $E$ and a strict \nbd{E_0}semigroup $\vt$ on $\sB^a(E)$ such that the product system of $\vt$ is $E^\odot$.
\ethm

\brem\label{genrem}
Note that both theorems remain true whenever for one member $E_t$ of $E^\odot$ with $t\ne0$, a suitable multiple $E_t^\en$ ($\en$ a cardinal number) of $E_t$ has a direct summand $\cB$, also if $\cB$ is not \nbd{\sigma}unital.
\erem

\section{The continuous case}\label{contSEC}

We now switch to the problem when, in the situation of Theorem \ref{algcthm}, the product system is also continuous. Following the same idea as described in \eqref{idea}, we simply could pass to the Lebesgue measure, show that the direct integrals make sense, and convince ourselves that the resulting \nbd{E_0}semigroup is strongly continuous and gives back the continuous structure on $E^\odot$ we started with. This is possible, but the operations mod $1$ create quite horrible problems to write it down.

In Skeide \cite{Ske07}, in the case of unital $\cB$, we followed a different idea due to Arveson \cite{Arv06} in the case of Hilbert spaces. Suppose $E_1$ has unit vector $\xi_1$. (For Hilbert spaces this is not a problem. For continuous full product systems and unital $\cB$ this is automatic; see \cite[Lemma 3.2]{Ske07}.) Consider those sections $x=\bfam{x_\alpha}_{\alpha\in\R_+}$ of $E^\odot$ that are \it{stable} with respect to $\xi_1$ in the sense that there exists $\alpha_0\in\R_+$ such that $\xi_1x_\alpha=x_{\alpha+1}$ for all $\alpha\ge\alpha_0$. Then one may define a semiinner product on the space of stable sections by setting $\AB{x,y}=\lim_{T\to\infty}\int_T^{T+1}\AB{x_\alpha,y_\alpha}\,d\alpha$. (This limit is over a function of $T$ which is eventually constant.) One may divide out the kernel of that inner product and complete. On that space  of sections mod $\AB{\bullet,\bullet}$ the product system ``acts'' from the right as $xy_t=\bfam{x_{\alpha-t}y_t}_{\alpha\in\R_+}$ (where we put $x_\alpha=0$ for negative $\alpha$).

This approach yields the same result, as if we constructed the dilation of the discrete subsystem based on the unit vector $\xi_1$ and used it as input for \eqref{idea}. (We proved that in Skeide \cite{Ske06a} for Hilbert spaces.) But it does not have anymore the problems with addition on $\RO{0,1}$ mod $1$. It is not too much trouble to prove the desired continuity results. We shall try to see how Arveson's approach can be safed for the nonunital case.

Suppose $E^\odot$ is a full continuous product system, and suppose (following Remark \ref{genrem}) that for some $t\ne0$ and for some cardinal number $\en$ the multiple $E_t^\en$ has a direct summand $\cB$. By rescaling, we may assume that $t=1$. Once more, the correspondence $\sB^a(\cB^\en,E_1^\en)$ over $\sB^a(\cB^\en)$ has a unit vector $\Xi_1$, say. This vector cannot act on $E_\alpha$. It can, however, act on $E_\alpha^\en$. Let $S$ be a set of cardinality $\#S=\en$ and denote the elements of $E_\alpha^\en$ as $X_\alpha=\bfam{X_\alpha^s}_{s\in S}$. Then put $\Xi_1X_\alpha:=\bfam{\sum_{s'\in S}(\Xi_1)_{ss'}X_\alpha^{s'}}_{s\in S}$.

We start by defining the direct integrals we need. Let the continuous structure of $E^\odot$ be determined by the family $i$ of embeddings $i_t\colon E_t\rightarrow\wh{E}$. This gives to embeddings $i_t^\en\colon E_t^\en\rightarrow\wh{E}^\en$. Every section $X=\bfam{X_t}_{t\in\R_+}$ with $X_t\in E_t^\en$ gives rise to a function $t\mapsto X(t):=i_t^\en X_t$ with values in $\wh{E}^\en$. We denote by
\beqn{
CS_i^\en(E^\odot)
~=~
\BCB{\,X\colon t\mapsto X(t)\text{~is continuous~}}
}\eeqn
the set of all such sections that are continuous. Let $0\le a<b<\infty$. By $\int_a^b E_\alpha^\en\,d\alpha$ we understand the norm completion of the pre-Hilbert \nbd{\cB}module that consists of continuous sections $X\in CS_i^\en(E^\odot)$ restricted to $\RO{a,b}$ with inner product
\beqn{
\AB{X,Y}_{\SB{a,b}}
~:=~
\int_a^b\AB{X_\alpha,Y_\alpha}\,d\alpha
~=~
\int_a^b\AB{X(\alpha),Y(\alpha)}\,d\alpha.
}\eeqn
Note that all continuous sections are bounded on the compact interval $\SB{a,b}$ and, therefore, square integrable. The following proposition is proved precisely as \cite[Proposition 4.2]{Ske07} (which holds for arbitrary subbundles of Banach bundles).

The following proposition is proved as \cite[Proposition 4.3]{Ske07}.	

\bprop\label{rcllprop}
$\int_a^b E_\alpha^\en\,d\alpha$ contains the space $\eR_{\RO{a,b}}$ of restrictions to $\RO{a,b}$ of those sections $X$ for which $t\mapsto X(t)$ is right continuous with finite jumps (this implies that there exists a left limit) in finitely many points and bounded on $\RO{a,b}$, as a pre-Hilbert submodule.
\eprop

Let $\sS$ denote the right \nbd{\cB}module of all sections $X$ that are \hl{locally $\eR$}, that is, for every $0\le a<b<\infty$ the restriction of $X$ to $\RO{a,b}$ is in $\eR_{\RO{a,b}}$, and which are \hl{stable} with respect to the unit vector $\Xi_1$ in $\sB^a(\cB^\en,E_1^\en)$, that is, there exists an $\alpha_0\ge0$ such that
\beqn{
X_{\alpha+1}
~=~
\Xi_1X_\alpha
}\eeqn
for all $\alpha\ge\alpha_0$. By $\sN$ we denote the subspace of all sections in $\sS$ which are eventually $0$, that is, of all sections $X\in\sS$ for which there exists an $\alpha_0\ge0$ such that $X_\alpha=0$ for all $\alpha\ge\alpha_0$. A straightforward verification shows that
\beqn{
\AB{X,Y}
~:=~
\lim_{m\to\infty}\int_m^{m+1}\AB{X(\alpha),Y(\alpha)}\,d\alpha
}\eeqn
defines a semiinner product on $\sS$ and that $\AB{X,X}=0$ if and only if $X\in\sN$. Actually, we have
\beqn{
\AB{X,Y}
~=~
\int_T^{T+1}\AB{X(\alpha),Y(\alpha)}\,d\alpha
}\eeqn
for all sufficiently large $T>0$; see \cite[Lemma 2.1]{Arv06}. So, $\sS/\sN$ becomes a pre-Hilbert module with inner product $\AB{X+\sN,Y+\sN}:=\AB{X,Y}$. By $E$ we denote its completion.

\bprop\label{denseprop}
For every section $X$ and every $\alpha_0\ge0$ define the section $X^{\alpha_0}$ as
\beqn{
X^{\alpha_0}_\alpha
~:=~
\begin{cases}
0&\alpha<\alpha_0
\\
\Xi_1^nX_{\alpha-n}&\alpha\in\RO{\alpha_0+n,\alpha_0+n+1},n\in\N_0.
\end{cases}
}\eeqn

If $X$ is in $CS_i^\en(E^\odot)$, then $X^{\alpha_0}$ is in $\sS$. Moreover, the set $\bCB{X^{\alpha_0}+\sN\colon X\in CS_i^\en(E^\odot),\alpha_0\ge0}$ is a dense submodule of $E$.
\eprop

After these preparations it is completely plain to see that for every $t\in\R_+$ the map $X\odot y_t\mapsto Xy_t$, where
\beqn{
(Xy_t)_\alpha
~=~
\begin{cases}
X_{\alpha-t}y_t&\alpha\ge t,
\\
0&\text{else},
\end{cases}
}\eeqn
and where $X_\alpha y_t=\bfam{X_\alpha^sy_t}{s\in S}$, defines an isometry $v_t\colon E\odot E_t\rightarrow E$, and that these isometries iterate associatively.

\bprop\label{surprop}
Each $v_t$ is surjective.
\eprop

\proof
By Proposition \ref{denseprop} it is sufficient to approximate every section of the form $X^{\alpha_0}$ with $X\in CS_i^\en(E^\odot),\alpha_0\ge0$ in the (semi-)inner product of $\sS$ by finite sums of sections of the form $Yz_t$ for $Y\in\sS,z_t\in E_t$. As what the section does on the finite interval $\RO{0,t}$ is not important for the inner product, we may even assume that $\alpha_0\ge t$. And as in the proof of Proposition \ref{denseprop} the approximation can be done by approximating $X$ in $\eR_{\RO{\alpha_0,\alpha_0+1}}$ and then extending the restriction to $\RO{\alpha_0,\alpha_0+1}$ stably to the whole axis. (This stable extension is the main reason why we worry to introduce the subspace of right continuous sections.)

Proposition \ref{surprop} for $\en=1$ and unital $\cB$ is done in \cite[Proposition 4.6]{Ske07}. The restriction that $\cB$ be unital can be omitted without affecting the proof. One may either repeat the proof word by word (for functions with values in $\wh{E}^\en$ instead of $\wh{E}$). Or one may note that the proof goes through for any finite $\en$, and that the approximation maybe done (with one more $\ve$) by restricting to a suitable finite subset $S$ (of course, depending on the section to be approximated).\qed

\lf
So, the $v_t$ form a dilation of $E^\odot$ to $E$. Like in \cite[Proposition 4.7]{Ske07}, we show that the dilation is \hl{continuous} in the following sense.

\bprop\label{cdprop}
For every $x\in E$ and every continuous section $y\in CS_i(E^\odot)$ the function $t\mapsto xy_t$ is continuous.
\eprop

In the proof of \cite[Proposition 4.7]{Ske07} just replace the section $z\in CS_i(E^\odot)$ by $Z\in CS_i^\en(E^\odot)$.

\bcor
The \nbd{E_0}semigroup $\vt^v$ is strictly continuous.
\ecor

\proof
This proof is almost identical to that of \cite[Corollary 4.8]{Ske07}. The only problem is that in our context here we do not have available a continuous section $\zeta$ of unit vectors that would fulfill $x\zeta_\vt\to x\U=x$ for all $x\in E$. Instead, for a given $x\in E$ we choose $\beta\in\cB$ with $\norm{\beta}\le1$ and $x\beta$ sufficiently close to $x$. Then we choose a continuous section $\zeta\in CS_i(E^\odot)$ with $\zeta_0=\beta$. With that section $\zeta$ everything goes exactly like in the proof of \cite[Corollary 4.8]{Ske07}.\qed

\bcor
The continuous structure induced by the \nbd{E_0}semigroup $\vt^v$ coincides with the continuous structure of $E^\odot$.
\ecor

\proof
This is Proposition \ref{cdprop} together with Theorem \ref{cisothm}.\qed

\brem
For unital $\cB$, this corollary is \cite[Proposition 4.9]{Ske07}. In the proof of \cite[Proposition 4.9]{Ske07} we showed, however, only right continuity. Theorem \ref{cisothm} settles this gap.
\erem

\brem
Note that the proofs of the two preceding corollaries do not depend on the concrete form of the left dilation. We, therefore, showed the following more general statement: If $v_t$ is a left dilation of a continuous product system $E^\odot$ that is continuous in the sense of Proposition \ref{cdprop}, then the induced \nbd{E_0}semigroup $\vt^v$ is strongly continuous and the continuous structure induced by that \nbd{E_0}semigroup coincides the original one.
\erem

We summarize.

\bthm\label{mthm}
Every full continuous product system of correspondences over a \nbd{\sigma}unital \nbd{C^*}al\-ge\-bra $\cB$ is the continuous product system associated with a strictly continuous \nbd{E_0}semigroup that acts on the algebra of all adjointable operators on a full Hilbert \nbd{\cB}module.
\ethm

For the sake of completeness, we state in the \nbd{\sigma}unital case the classification theorem for \nbd{E_0}semigroup by product systems as stated in Skeide \cite{Ske08p1} for unital \nbd{C^*}algebras. We refer to \cite{Ske08p1} for the definition of stable cocycle conjugacy. Precisely under the same conditions as in \cite[Section 9]{Ske08p1} for unital \nbd{C^*}algebras $\cB$, we obtain the following theorem for \nbd{\sigma}unital \nbd{C^*}algebras. Recall that a continuous product system $E^\odot$ is \hl{countably generated}, if there is a countable subset $S$ of $CS_i(E^\odot)$ such that $CS_i(E^\odot)$ is the locally uniform closure of the linear span of $S$.

\bthm\label{C*contclassthm}
Let $\cB$ be a \nbd{\sigma}unital \nbd{C^*}algebra. Then there is a one-to-one correspondence between equivalence classes (up to stable cocycle conjugacy with strongly continuous cocycles) of strongly continuous strict \nbd{E_0}semigroups acting on the operators of countably generated full Hilbert \nbd{\cB}mod\-ules and isomorphism classes of countably generated continuous product systems of full correspondences over $\cB$.
\ethm


\setlength{\baselineskip}{2.5ex}


\newcommand{\Swap}[2]{#2#1}\newcommand{\Sort}[1]{}
\providecommand{\bysame}{\leavevmode\hbox to3em{\hrulefill}\thinspace}
\providecommand{\MR}{\relax\ifhmode\unskip\space\fi MR }
\providecommand{\MRhref}[2]{%
  \href{http://www.ams.org/mathscinet-getitem?mr=#1}{#2}
}
\providecommand{\href}[2]{#2}


\end{document}